\documentclass[reqno,12pt]{amsart}

\usepackage{amsthm,amsfonts}

\usepackage[centertags]{amsmath}

\usepackage[a4paper,margin=34mm]{geometry}

\usepackage{amssymb,latexsym}

\usepackage[centertags]{amsmath}

\usepackage{bm}

\newtheorem{theorem}{Theorem}
\newtheorem{proposition}{Proposition}

\newtheorem{definition}{Definition}

\newtheorem{remark}{Remark}

\begin{document}

\title[Full Randomness in Higher Difference Structure of Markov Chains]{Full Randomness in the Higher Difference Structure of Two-state Markov Chains}

\author[A. Yu. Shahverdian]{A. Yu. Shahverdian}

\address{Institute for Informatics and Automation Problems of NAS RA}

\email{svrdn@yerphi.am}



\date{\today}

\maketitle

\begin{abstract}
The paper studies the higher-order absolute differences taken from progressive terms of time-homogenous binary Markov chains. Two theorems presented are the limiting theorems for these differences, when their order $k$ converges to infinity. Theorems 1 and 2 assert that there exist some infinite subsets $E$ of natural series such that $k$th order differences of every such chain converge to the equi-distributed random binary process as $k$ growth to infinity remaining on $E$. The chains are classified into two types and $E$ depend only on the type of a given chain. Two kinds of discrete capacities for subsets of natural series are defined, and in their terms such sets $E$ are described.
\end{abstract}

\section{\textbf{Introduction}}\label{sec1}

In this paper an application of the suggested in \cite{1}-\cite{8} difference analysis to studying binary Markov chains is presented. In difference analysis we are interested in the following question: which is the higher-order difference structure of a given process and how this structure can characterize the process.

The paper studies the time-homogenous binary Markov chains $\bm{\xi} =(\xi_{n})_{n\geq 0}$ where every $\xi_{n}$ is binary variable describing the state of the cain $\bm{\xi}$ at the moment $n$. The main results, Theorem~\ref{T1} and Theorem \ref{T2}, are the limiting theorems for such chains. These theorems concern infinite sets $E\subseteq \mathbb{N}$ which posses such a property: for arbitrary chain $\bm{\xi}$ the $E$  permits the existence of the limit of $k$th order absolute differences $\bm{\xi}^{(k)}$, when $k$ converges to $\infty$ remaining on $E$. The existence of such $E$ is claimed and their description in capacity terms is given. The chains are classified into two types and $E$ depend only on the type of a given chain. The limiting process is the equi-distributed random binary sequence, denoted $\bm{\theta}$ (see Eq.~(\ref{E2})).

Theorems \ref{T1} and \ref{T2} improve our previous results from \cite{7,8}; some details on this matter in Section \ref{sec3} (points {\em (a) - (c)}) are given. The limiting process, which is the equi-distributed sequence, should be recognized as the most random binary sequence. Therefore, theorems presented state the existence of full randomness in the higher difference structure of arbitrary time-homogenous binary Markov chain.

Let us explain our statement in more detail. Let
$$
\bm{\xi} = (\xi_{0}, \xi_{1}, \ldots, \xi_{n}, \ldots)
$$
be some random sequence whose components $\xi_{n}$ take binary values $x\in X$, $X = \{0, 1\}$
with some positive probabilities, $P(\xi_{n}= x) = p_{n}(x)$. Then $k$th order ($k\geq 0$) absolute differences $\xi_{n}^{(k)}$, which are defined recurrently,
$$
\xi_{n}^{(0)} \equiv \xi_{n} \quad \text{and} \quad \xi_{n}^{(k)} = |\xi_{n+1}^{(k-1)} - \xi_{n}^{(k-1)}|
\quad (n\geq 0),
$$
also take binary values with some probabilities $P(\xi_{n}^{(k)}= x) = p_{n}^{(k)}(x)$,
and hence, one can consider $k$th order difference random binary sequence
$$
\bm{\xi}^{(k)} = (\xi_{0}^{(k)}, \xi_{1}^{(k)}, \ldots, \xi_{n}^{(k)}, \ldots ).
$$
Our interest is the limits of $\bm{\xi}^{(k)}$ when $k$ goes to infinity. Let some infinite $E \subseteq \mathbb{N}$ be given. We say that $\bm{\xi}^{(k)}$ converge on $E$ to a random binary sequence $\bm{\xi}_{E}^{\infty}$,
and denote this
$$
\bm{\xi}^{\infty}_{E} = \lim_{\substack{k\to \infty \\ k\in E}}\bm{\xi}^{(k)},
$$
if for $n\in \mathbb{N}$ and $x\in X$ the probabilities $p_{n}^{(k)}(x)$  tend to some numbers $p_{n}^{(\infty)}(x)$ as $k\to \infty$ and $k\in E$,
$$
\lim_{\substack{k\to \infty \\ k\in E}}p_{n}^{(k)}(x)  = p^{(\infty)}_{n}(x)
$$
({\em convergence by probability} on $E$ and {\em partial limits}). Therefore, $\bm{\xi}^{\infty}_{E}$ is a random binary sequence,
$$
\bm{\xi}_{E}^{\infty} = (\xi_{0}^{(\infty)}, \xi_{1}^{(\infty)}, \ldots, \xi_{n}^{(\infty)}, \ldots)
$$
whose components $\xi_{n}^{(\infty)}$ take the values $x\in X$ with probabilities $P(\xi_{n}^{(\infty)}= x) = p_{n}^{(\infty)}(x)$ (which depend on $E$).

We consider binary Markov chains $\boldsymbol{\xi}= (\xi_{n})_{n=0}^{\infty}$ whose state space $X$ consists of two binary symbols, $X =\{0,1\}$. We assume that the chains $\bm{\xi}$ are time-homogeneous, that is, for $x$, $x_i$, $y \in X$
\begin{align}\label{e1}
&P(\xi_{n}=y|\xi_{n-1}=x, \xi_{n-2}=x_{1}, \ldots, \xi_{0}=x_{n-1}) \nonumber \\
&= P(\xi_{n}=y|\xi_{n-1}=x) \nonumber
\end{align}
(Markov property) and there is some function $\pi(x,y)$ on $X\times X$ such that for  $n\geq 1$ and $x,y \in X$
$$
P(\xi_{n}=y|\xi_{n-1}= x) = \pi(x,y)
$$
(homogeneity: one-step transition probabilities $P(\xi_{n}=y|\xi_{n-1}=x)$ do not depend on time $n$). It is also assumed that some initial distribution of probabilities $P(\xi_{0} =x)$ on $X$ is given. In what follows it is always assumed that $\bm{\xi}$ denotes the time-homogeneous binary Markov chain.

Some simple computations testify, that if for given $\bm{\xi}$ an infinite $E\subseteq \mathbb{N}$ is chosen arbitrarily, then the limiting process $\bm{\xi}^{\infty}_{E}$ may not exist. On the other hand, it follows from \cite{7,8} that for $E = \{2^{m}-1: m\geq 0\}$ and large collection of $\bm{\xi}$ the limit $\bm{\xi}^{\infty}_{E}$ exists and it is equi-distributed random sequence. The problem studied relates to the following question: how the sets $E\subseteq \mathbb{N}$, which for arbitrary Markov chain $\bm{\xi}$ permit the existence of $\bm{\xi}^{\infty}_{E}$, can be described?

The main results of this paper, Theorems~\ref{T1} and \ref{T2}, are the limiting theorems for such chains.
Two discrete capacities for subsets of $\mathbb{N}$ are defined and in their terms such sets $E$ are described.
The limiting process $\bm{\xi}_{E}^{\infty}$, whose existence assert these theorems, is the equidistributed
random binary sequence.

The (stochastic) transition matrix of every time homogenous binary Markov chain $\bm{\xi}$ can be written as
$$
Q_{\bm{\xi}} =
\begin{bmatrix}
s      &  1-s  \\
1-p    &  p    \\
\end{bmatrix}
\qquad (0<s,p<1)
$$
where $s = \pi(0,0)$ and $p= \pi(1,1)$. Theorems \ref{T1} and \ref{T2} consider two types of chains $\bm{\xi}$, depending on which of the next two relationships (I) and (II) between $s$ and $p$
\begin{equation}\label{E1}
\text{(\textbf{I})} \ \  s\neq p \ \ \text{and} \ \ s\neq 1-p,  \ \ \text{(\textbf{II})} \ \  s=p \ \ \text{or} \ \ s=1-p
\end{equation}
holds: we say that the chain $\bm{\xi}$ is of I-st or II-nd type whenever for $s$ and $p$ the relations (I) or (II) (respectively) from Eq.~(\ref{E1}) are satisfied.

The paper consists of four sections. The next Section~\ref{sec2} contains definitions of discrete capacities that we use. In Section~\ref{sec3} the formulations of main Theorems \ref{T1} and \ref{T2} are presented, and last Section~\ref{sec4} contains some additional comments.

\section{\textbf{Some definitions}}\label{sec2}

To proceed to formulation of our Theorems~\ref{T1} and \ref{T2}, we need to present two discrete capacities  $\mathcal{C}$ and $\bm{c}$ defined for subsets of natural series $\mathbb{N}$. Their definition is given by means of binary representation of natural numbers and binary version of Pascal triangle.
The binary Pascal triangle $\mathbb{P}$ and its $k$th line $\bm{\ell}_k$ are defined as
$$
\mathbb{P}= \{\alpha_{k,i}: k\geq 0, 0\leq i \leq k\}, \ \bm{\ell}_{k} = (\alpha_{k,0},  \alpha_{k,1}, \ldots, \alpha_{k,k})
$$
that is, $ \mathbb{P} = \bigcup_{k=0}^{\infty} \bm{\ell}_{k}$; here, $\alpha_{k,i}\in \{0,1\}$ are the following: $\alpha_{0,0} = 1$ (the vertex of $\mathbb{P}$ and the line $\bm{\ell}_{0}$), $\alpha_{1,0} = \alpha_{1,1} =1$ (the line $\bm{\ell}_{1}$), and for $k\geq 2$ the line $\bm{\ell}_{k}$ consists of such $\alpha_{k,i}$,
$$
\alpha_{k,i} = \left \{
\begin{array}{ll}
0, & \binom{k}{i} \ \text{is even} \\[0.25cm]
1, & \binom{k}{i} \ \text{is odd}
\end{array}
\qquad (0\leq i\leq k).
\right.
$$
One can see that this is the same as if for $k\geq 1$ one defines:
$\alpha_{k,0} = \alpha_{k,k}=1$, and
$$
\alpha_{k,i} = |\alpha_{k-1,i-1} - \alpha_{k-1,i}| \qquad (1\leq i\leq k-1).
$$

The capacities $\mathcal{C}$ and $\bm{c}$ are defined by means of some quantities related to binary expansion of natural numbers. For $k\geq 1$ its binary representation is given as
\begin{equation}\label{E1a}
k = \sum_{i=0}^{p}\varepsilon_{i}2^{i} \quad \text{where} \ \ p\geq 0, \ \ \varepsilon_{i}\in \{0,1\},
\ \  \varepsilon_{p} =1;
\end{equation}
we denote
$$
b(k) = \sum_{i=0}^{p}\varepsilon_{i},
\quad
\beta(k) = \sum_{i=0}^{k}\alpha_{k,i}.
$$
For natural $k$ we use the following notations: $\nu(k)$ denotes the maximal of such $m$, $0\leq m \leq p$ for which all the coefficients $\varepsilon_{i}$, $0\leq i\leq m$ in expansion (\ref{E1a}) are equal to $1$,
$$
\nu(k) = \max\{m: \varepsilon_{0} = \varepsilon_{1} = \cdots =\varepsilon_{m} =1\};
$$
$\mu(k)$ denotes the maximal of such $m$, $0\leq m \leq k$ for which all the $\alpha_{k,i}$, $0\leq i\leq m$ (first $m$ entries of the line $\bm{\ell}_{k}$ of the triangle $\mathbb{P}$) are equal $1$,
$$
\mu(k) = \max\{m: \alpha_{k,0} = \alpha_{k,1} = \cdots = \alpha_{k,m} =1\}.
$$
The capacities $\mathcal{C}$ and $\bm{c}$ are assigned on the collection $2^{\mathbb{N}}$ of subsets of natural series and defined as follows.
\begin{definition}
For $e\subseteq \mathbb{N}$ we define
$$
\mathcal{C}(e) = \sum_{k \in e}\nu(k), \quad \bm{c}(e) = \sum_{k \in e}b(k).
$$
\end{definition}

The $\mathcal{C}(e)$ and $\bm{c}(e)$ can be expressed by the entries of the Pascal triangle $\mathbb{P}$:
one can prove that $\mu(k) = 2^{\nu(k)}$ and $\beta(k) = 2^{b(k)}$, and, therefore,
$$
 \mathcal{C}(e) = \sum_{k\in e}\log_{2}\mu(k),  \quad \bm{c}(e) = \sum_{k\in e}\log_{2}\beta(k).
$$
Both $\mathcal{C}$ and $\bm{c}$ are differed from discrete capacity, considered in denumerable Markov chains and random walk (e.g., \cite{9}); for details on $\mathcal{C}$ and $\bm{c}$ see \cite{4} and \cite{8}. We denote $\mathcal{C}(k)= \mathcal{C}(\{k\})$ and $\bm{c}(k) = \bm{c}(\{k\})$.

Let us present an example of computation of these capacities. We denote $L_{p} = \{k\in \mathbb{N}: 2^{p-1} \leq k < 2^{p}\}$ and for $p\geq 2$ and $0\leq s \leq p$ consider the sets $B_{p}(s)$ and $b_{p}(s)$:
$$
B_{p}(s) = \{k\in L_{p}: \nu_{k}\geq s\}, \ \  b_{p}(s)= \{k\in L_{p}: b(k)\geq s\}.
$$
The complement of $b_{p}(s)$  is the Hamming ball of radius $s$  (\cite{4}; there is a misprint in \cite{4} on computation of capacity of these balls).

\begin{proposition}\label{PP2}
For $p\geq 2$ and $0\leq s \leq p$ the relations
$$
\mathcal{C}(B_{p}(s)) = \sum_{i=0}^{s}i2^{p-i}, \quad \bm{c}(b_{p}(s)) = \sum_{i=s}^{p}i\binom{p}{i}
$$
are true.
\end{proposition}

\section{\textbf{Main theorems}}\label{sec3}

In this section we formulate our main results, Theorems~\ref{T1} and \ref{T2}. They define some sets $E\subseteq \mathbb{N}$ and state the convergence of $k$th order difference processes $\boldsymbol{\xi}^{(k)}$ (as $k\to \infty$ and $k\in E$) to the equi-distributed random binary sequence $\bm{\theta}$; the $\bm{\theta}$ is defined as
\begin{equation}\label{E2}
\bm{\theta} = (\theta_{0}, \theta_{1}, \ldots, \theta_{n}, \ldots) \quad \text{where} \quad  P(\theta_{n}=x) = \frac{1}{2}
\end{equation}
for all $n\geq 0$ and $x\in \{0,1\}$. In next formulations $o_{k}(1)$ denotes the Landau symbol: it is some numerical quantity which tends to $0$ as $k$ converges to $\infty$.

Theorem~\ref{T1} and Theorem~\ref{T2}, formulated in next subsections, improve some our results from \cite{7,8}. If compared with \cite{7,8}, the improvement is due to the following three features of Theorems~\ref{T1} and \ref{T2}: \ {\em (a)} \ the sets $E$ in formulations of these theorems depend only on the type (I-st or II-nd type) of the chain $\bm{\xi}$ and do not depend on other details uniquely determining a given chain; \ {\em (b)} \ the theorems estimate the rate (exponential) of the convergence; \ {\em (c)} \ a different description of sets $E$ (Eqs.~(\ref{E3}) and (\ref{E6})) is given.

In next Sections~\ref{sec3.1} and \ref{sec3.2} we present some examples of such sets $E$ (Eqs.~(\ref{E5}) and (\ref{E8})). These examples appear to be quite general and connect us (Propositions~\ref{P2} and \ref{P4}) with another, considered in \cite{8}, description of these sets. In addition, Remarks~\ref{R1} and \ref{R2} state that the sets $E$ from these examples are the 'largest' ones, satisfying the assumptions (\ref{E3}) and (\ref{E6}) in these theorems. This allows  us to derive some conclusions (Propositions~\ref{P3} and \ref{P5}) on densities of sets $E$ from Theorems~\ref{T1} and \ref{T2}.

\subsection{Chains of I-st type}\label{sec3.1}

Let us formulate our Theorem \ref{T1} which concerns Markov chains of I-st type (defined by Eq.~(\ref{E1})). This theorem describes infinite sets $E\subseteq \mathbb{N}$ which possess the property that the limiting processes $\boldsymbol{\xi}^{\infty}_{E}$ exists for arbitrary Markov chain $\boldsymbol{\xi}$ of I-st type: the theorem asserts the convergence of $k$th order difference processes $\boldsymbol{\xi}^{(k)}$ (as $k\to \infty$ and $k\in E$) to the equi-distributed process $\bm{\theta}$ (defined by Eq.~(\ref{E2})).

\begin{theorem}\label{T1}
Let a set $E\subseteq \mathbb{N}$ be such that
\begin{equation}\label{E3}
\lim_{\substack{k\to \infty \\ k\in E}}\mathcal{C}(k) = \infty
\end{equation}
and $\bm{\xi}$ be Markov chain of I-st type. Then the limiting process $\bm{\xi}^{\infty}_{E}$ exists and
$\bm{\xi}^{\infty}_{E}= \bm{\theta}$, that is,
\begin{equation}\label{E4}
\lim_{\substack{k \to \infty \\ k\in E}} \bm{\xi}^{(k)} = \bm{\theta}.
\end{equation}
The convergence in Eq.~(\ref{E4}) is exponential: given $\bm{\xi}$ there is some $\delta$, $|\delta| <1$ which depends only on transition matrix of $\bm \xi$, such that for $n\geq 1$, $k\in E$ and $\lambda \in \{0,1\}$ the relation
$$
P(\xi_{n}^{(k)} =\lambda) = \frac{1}{2} + o_{k}(1)\delta^{k}
$$
holds.
\end{theorem}

Let us present some examples of sets $E\subseteq \mathbb{N}$ satisfying Eq.~(\ref{E3}). With this aim we consider the unions of $B_{p}(s_{p})$,
\begin{equation}\label{E5}
E = \bigcup_{p=1}^{\infty}B_{p}(s_{p}).
\end{equation}

\begin{proposition}\label{P2}
Let $E \subseteq \mathbb{N}$ be defined by Eq.~(\ref{E5}). Then $E$ satisfies  Eq.~(\ref{E3})
if and only if the conditions
\begin{equation}\label{E5a}
\lim_{p\to \infty}s_{p} =\infty \quad \text{and} \quad  \sum_{p=1}^{\infty}2^{-p}\mathcal{C}(B_{p}(s_{p})) = \infty
\end{equation}
hold.
\end{proposition}

The next Remark~\ref{R1} asserts that the given by Proposition~\ref{P2} example of sets $E$ satisfying Eq.~(\ref{E3}) is quite general.

\begin{remark}\label{R1}
For a set $E\subseteq \mathbb{N}$ the condition (\ref{E3}) holds if and only if there is a set $E^{\prime}\subseteq \mathbb{N}$ of the form (\ref{E5}) satisfying (\ref{E3}) and such that $E\subseteq E^{\prime}$.
\end{remark}

We describe the density of sets $E$ from Theorem~\ref{T1}. For $m\geq 1$ we denote $E_{m} = \{k\in E: 1\leq k\leq m\}$, consider the ratio $\rho_{m}(E) = \dfrac{|E_{m}|}{m}$,  where $|E_{m}|$ denotes the cardinality of $E_{m}$, and define
$$
dens(E) = \lim_{m\to \infty}\rho_{m}(E).
$$

\begin{proposition}\label{P3}
If a set $E\subseteq \mathbb{N}$ satisfies Eq.~(\ref{E3}), then $dens(E) =0$. For a given $0< \delta_{m} \leq 1$, $\delta_{m} \downarrow 0$ there is a set $E \subseteq \mathbb{N}$ which satisfies Eq.~(\ref{E3}) and such that
$\rho_{m}(E) \geq \delta_{m}$ for all $m\geq 1$.
\end{proposition}

\subsection{Chains of II-nd type}\label{sec3.2}

The next Theorem~\ref{T2} concerns Markov chains of II-nd type and describes infinite sets $E\subseteq \mathbb{N}$, which possess the property that the limiting processes $\boldsymbol{\xi}^{\infty}_{E}$ exists for arbitrary Markov chains $\boldsymbol{\xi}$ of II-nd type; the limiting process is again the equi-distributed process $\bm{\theta}$.

\begin{theorem}\label{T2}
Let a set $E\subseteq \mathbb{N}$ be such that
\begin{equation}\label{E6}
\lim_{\substack{k\to \infty \\ k\in E}}\bm{c}(k) = \infty
\end{equation}
and $\bm{\xi}$ be Markov chain of II-nd type. Then the limiting process $\bm{\xi}^{\infty}_{E}$ exists and
$\bm{\xi}^{\infty}_{E}= \bm{\theta}$, that is,
\begin{equation}\label{E7}
\lim_{\substack{k \to \infty \\ k\in E}} \bm{\xi}^{(k)} = \bm{\theta}.
\end{equation}
The convergence in Eq.~(\ref{E7}) is exponential: given $\bm{\xi}$ there is some $\delta$, $|\delta| <1$ which depends only on transition matrix of $\bm \xi$, such that for $n\geq 1$, $k\in E$ and $\lambda \in \{0,1\}$ the relation
$$
P(\xi_{n}^{(k)} =\lambda) = \frac{1}{2} + o_{k}(1)\delta^{k}
$$
holds.
\end{theorem}

As the examples of sets $E\subset \mathbb{N}$ satisfying Eq.~(\ref{E6}) we consider the unions of
$b_{p}(s_{p})$,
\begin{equation}\label{E8}
E = \bigcup_{p=1}^{\infty}b_{p}(s_{p}).
\end{equation}

\begin{proposition}\label{P4}
Let $E\subseteq \mathbb{N}$ be defined by Eq.~(\ref{E8}). Then $E$ satisfies Eq.~(\ref{E6})
if and only if the conditions
\begin{equation}\label{E10}
\lim_{p\to \infty}s_{p} =\infty \quad \text{and} \quad  \sum_{p=1}^{\infty}2^{-p}\bm{c}(b_{p}(s_{p})) = \infty
\end{equation}
hold.
\end{proposition}

\begin{remark}\label{R2}
For a set $E\subseteq \mathbb{N}$ the condition (\ref{E6}) holds if and only if there is a set $E^{\prime}\subseteq \mathbb{N}$ of the form (\ref{E8}) satisfying (\ref{E6}) and such that $E\subseteq E^{\prime}$.
\end{remark}

We compute the density of sets $E$ from Eq.~(\ref{E10}):
\begin{proposition}\label{P5}
If a set $E\subseteq \mathbb{N}$ defined by Eq.~(\ref{E8}) satisfies Eq.~(\ref{E10}), then $dens(E) =1$.
\end{proposition}

Particularly, Propositions~\ref{P4} and \ref{P5} imply that the sets $E$ from  Theorem~\ref{T2} can be as 'large', that their density equals $1$.

\section{\textbf{Some comments}}\label{sec4}

The chains considered can be treated as two-state probabilistic automata.

In \cite{4} independent random sequences have been studied (there is an unnecessary (and wrong) assumption in \cite{4} on independence of  $\bm{\xi}^{(k)}$). Theorem~\ref{T2} remains valid also for arbitrary independent identically distributed binary sequences.

The capacities $\mathcal{C}$ and $\bm{c}$ are some instances of the Fuglede-Choquet discrete capacities \cite{10, 11}, which are abstract version of classical capacities (e.g., \cite{11a}). Another kind of discrete capacity, applied to a  self-organized criticality model \cite{12}, is considered in \cite{13}.

The second relations in (\ref{E5a}) and (\ref{E10}) are the analogs for the Wiener criterion from potential theory (e.g., \cite{13a,14,15}; the sets $E$ from (\ref{E5}) and (\ref{E8}), satisfying these relations, are called thick sets. Apparently, the most known application of thick sets in classical theory is given by Keldysh theorem on the Dirichlet problem \cite{15}.

One of the basic concepts of ergodic theory is the notion of shift in probabilistic spaces \cite{16}.
E.g., independent sequences and Markov chains can be treated as consecutive iterates of some ergodic shifts (Bernoulli and Markov shifts \cite{16}). In \cite{4} we have defined the difference shift $M$; it is such, that
$k$th order absolute difference $\bm{\xi}^{(k)}$ coincides with $k$th iterate $M^{k}$ of the shift $M$ applied to the random sequence $\bm{\xi}$, $\bm{\xi}^{(k)}= M^{k}\bm{\xi}$. Therefore, Theorems~\ref{T1} and \ref{T2} can also be treated as some statements on iterates of the difference shift $M$.

\end{document}